\documentclass[12pt,italian]{article}

\title{\bf Probabilit\`a e paradosso}

\author{Germano D'Abramo$^1$ e Barbara D'Abramo$^2$\vspace{0.3cm}\\
{\small Istituto di Astrofisica Spaziale e Fisica Cosmica,}\\
{\small Roma, Italy}\\
{\small E--mail: {\tt dabramo@rm.iasf.cnr.it}}\vspace{0.3cm}\\
{\small $^2$Liceo Scientifico ``Leonardo Da Vinci'',}\\
{\small Fermo, Ascoli Piceno, Italy}\\
{\small E--mail: {\tt baxina@libero.it}}}

\date{\small 24 Giugno 2004}

\begin{document}

%

\maketitle

\section{Introduzione}

Come ognuno di voi ha avuto modo di sperimentare di persona, contare
manualmente una grande quantit\`a di oggetti \`e un'operazione faticosa e
per nulla esente da errori. La possibilit\`a di sbagliare \`e sempre in
agguato e la probabilit\`a di compiere un errore \`e significativamente
non nulla per gruppi di oggetti piuttosto numerosi.  Nella sezione~{\bf 2}
di questo breve articolo, giocando con le elementari regole del calcolo
delle probabilit\`a, mostriamo alcune sorprendenti e contro-intuitive
caratteristiche del processo di conteggio manuale. Infine, poniamo
all'attenzione del lettore una {\em rara} forma di inconoscibilit\`a
classica (cio\`e non quantistica), che ha qualcosa a che vedere con la
persecuzione della Sfortuna.

Nella sezione~{\bf 3}, invece, cerchiamo di dare una caratterizzazione
{\em statistica} del concetto di individuo `grande' o `geniale'. Cio\`e,
quale dovrebbe essere il numero medio di individui `grandi', rispetto al
numero medio che si otterrebbe nell'ipotesi che i grandi risultati si
verifichino {\em casualmente}, affinch\'e si possa parlare effettivamente
di `genio' o contributo `geniale'? Questa lettura del concetto di
genialit\`a, in stretta relazione con la nozione comune che si ha di esso,
sembra produrre una sorta di paradosso, che potrebbe essere definito come
il `paradosso del genio'.

Nella sezione~{\bf 4}, infine, descriviamo una semplice societ\`a ideale
nella quale viene messa a votazione democratica una perversa riforma degli
stipendi. Si vede, in questo caso, come la decisione pi\`u razionale da
prendere nell'esprimere il proprio voto democratico porti inevitabilmente
ad una conclusione perversa.

\section{Contare quante volte?}

Al compimento del quindicesimo anno, Sara ottiene dai suoi genitori il
permesso di rompere il maialino salvadanaio, da qualche tempo sofferente
di una forma (fatale) di costipazione da Euro.

Una volta compiuta l'inevitabile eutanasia, Sara tenta di contare il suo
cospicuo gruzzoletto, tutto costituito da monete di 1 Euro. Le conviene
contarlo una sola volta e fidarsi della cifra ottenuta o contarlo una
seconda volta per essere sicura del valore appena calcolato? La risposta
che viene spontanea un p\`o a tutti \`e la seconda, soprattutto se nel
conteggio si \`e avuta l'impressione di essersi distratti. Tuttavia,
vedremo ora che dietro questa scelta si nasconde una piccola sorpresa.

Supponiamo per semplicit\`a che la probabilit\`a che Sara ha di sbagliare
un conteggio sia $p$ e che sia sempre la stessa per tutti i conteggi
consecutivi che decide di fare. Quindi, se Sara effettua un solo conteggio
e poi decide di fidarsi del risultato, la probabilit\`a di aver sbagliato
\`e banalmente $P_{err1}=p$.

Se invece decidesse di ripetere una seconda volta il conteggio del suo
denaro, quale sarebbe la probabilit\`a di aver sbagliato il risultato
finale {\em almeno una volta}? \`E sufficiente un rapido calcolo per
convincersi che \`e $P_{err2}=1-(1-p)^2$, infatti $P_{err2}$ \`e uguale a
1 meno la probabilit\`a di {\em non} sbagliare in nessuno dei due
conteggi, cio\`e $(1-p)\cdot (1-p)$.

\begin{table}[h]
\begin{center}
\caption{Probabilit\`a di errore nei due casi:}
\begin{tabular}{ll}
 & \\
Un conteggio & $P_{err1}=p$\\
Due conteggi & $P_{err2}=1-(1-p)^2$
\end{tabular}
\end{center}
\end{table}

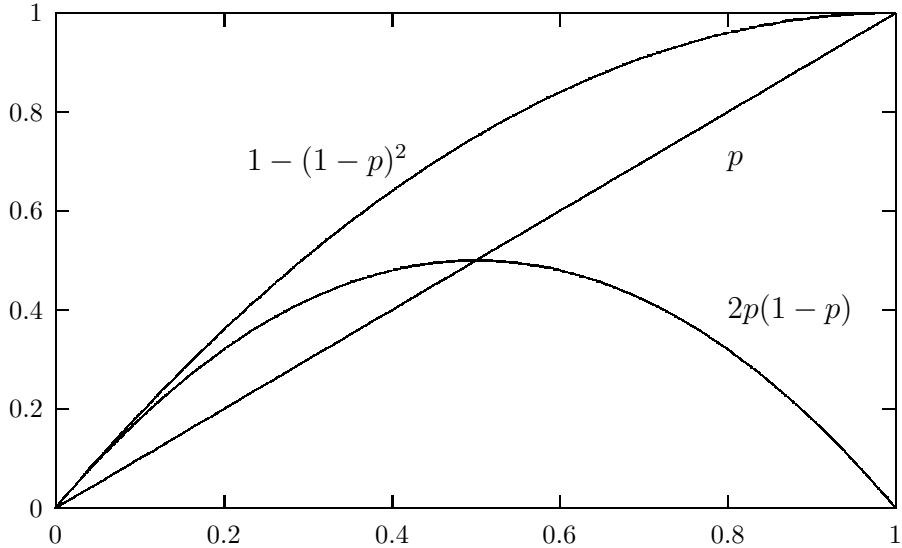
\begin{figure}

\setlength{\unitlength}{0.240900pt}
\ifx\plotpoint\undefined\newsavebox{\plotpoint}\fi
\sbox{\plotpoint}{\rule[-0.200pt]{0.400pt}{0.400pt}}%
\begin{picture}(1500,900)(0,0)
\font\gnuplot=cmr10 at 10pt
\gnuplot
\sbox{\plotpoint}{\rule[-0.200pt]{0.400pt}{0.400pt}}%
\put(120.0,82.0){\rule[-0.200pt]{4.818pt}{0.400pt}}
\put(100,82){\makebox(0,0)[r]{0}}
\put(1419.0,82.0){\rule[-0.200pt]{4.818pt}{0.400pt}}
\put(120.0,238.0){\rule[-0.200pt]{4.818pt}{0.400pt}}
\put(100,238){\makebox(0,0)[r]{0.2}}
\put(1419.0,238.0){\rule[-0.200pt]{4.818pt}{0.400pt}}
\put(120.0,393.0){\rule[-0.200pt]{4.818pt}{0.400pt}}
\put(100,393){\makebox(0,0)[r]{0.4}}
\put(1419.0,393.0){\rule[-0.200pt]{4.818pt}{0.400pt}}
\put(120.0,549.0){\rule[-0.200pt]{4.818pt}{0.400pt}}
\put(100,549){\makebox(0,0)[r]{0.6}}
\put(1419.0,549.0){\rule[-0.200pt]{4.818pt}{0.400pt}}
\put(120.0,704.0){\rule[-0.200pt]{4.818pt}{0.400pt}}
\put(100,704){\makebox(0,0)[r]{0.8}}
\put(1419.0,704.0){\rule[-0.200pt]{4.818pt}{0.400pt}}
\put(120.0,860.0){\rule[-0.200pt]{4.818pt}{0.400pt}}
\put(100,860){\makebox(0,0)[r]{1}}
\put(1419.0,860.0){\rule[-0.200pt]{4.818pt}{0.400pt}}
\put(120.0,82.0){\rule[-0.200pt]{0.400pt}{4.818pt}}
\put(120,41){\makebox(0,0){0}}
\put(120.0,840.0){\rule[-0.200pt]{0.400pt}{4.818pt}}
\put(384.0,82.0){\rule[-0.200pt]{0.400pt}{4.818pt}}
\put(384,41){\makebox(0,0){0.2}}
\put(384.0,840.0){\rule[-0.200pt]{0.400pt}{4.818pt}}
\put(648.0,82.0){\rule[-0.200pt]{0.400pt}{4.818pt}}
\put(648,41){\makebox(0,0){0.4}}
\put(648.0,840.0){\rule[-0.200pt]{0.400pt}{4.818pt}}
\put(911.0,82.0){\rule[-0.200pt]{0.400pt}{4.818pt}}
\put(911,41){\makebox(0,0){0.6}}
\put(911.0,840.0){\rule[-0.200pt]{0.400pt}{4.818pt}}
\put(1175.0,82.0){\rule[-0.200pt]{0.400pt}{4.818pt}}
\put(1175,41){\makebox(0,0){0.8}}
\put(1175.0,840.0){\rule[-0.200pt]{0.400pt}{4.818pt}}
\put(1439.0,82.0){\rule[-0.200pt]{0.400pt}{4.818pt}}
\put(1439,41){\makebox(0,0){1}}
\put(1439.0,840.0){\rule[-0.200pt]{0.400pt}{4.818pt}}
\put(120.0,82.0){\rule[-0.200pt]{317.747pt}{0.400pt}}
\put(1439.0,82.0){\rule[-0.200pt]{0.400pt}{187.420pt}}
\put(120.0,860.0){\rule[-0.200pt]{317.747pt}{0.400pt}}
\put(1175,627){\makebox(0,0)[l]{$p$}}
\put(420,627){\makebox(0,0)[l]{$1-(1-p)^2$}}
\put(1175,393){\makebox(0,0)[l]{$2p(1-p)$}}
\put(120.0,82.0){\rule[-0.200pt]{0.400pt}{187.420pt}}
\put(120,82){\usebox{\plotpoint}}
\multiput(120.00,82.59)(0.824,0.488){13}{\rule{0.750pt}{0.117pt}}
\multiput(120.00,81.17)(11.443,8.000){2}{\rule{0.375pt}{0.400pt}}
\multiput(133.00,90.59)(0.890,0.488){13}{\rule{0.800pt}{0.117pt}}
\multiput(133.00,89.17)(12.340,8.000){2}{\rule{0.400pt}{0.400pt}}
\multiput(147.00,98.59)(0.824,0.488){13}{\rule{0.750pt}{0.117pt}}
\multiput(147.00,97.17)(11.443,8.000){2}{\rule{0.375pt}{0.400pt}}
\multiput(160.00,106.59)(0.950,0.485){11}{\rule{0.843pt}{0.117pt}}
\multiput(160.00,105.17)(11.251,7.000){2}{\rule{0.421pt}{0.400pt}}
\multiput(173.00,113.59)(0.890,0.488){13}{\rule{0.800pt}{0.117pt}}
\multiput(173.00,112.17)(12.340,8.000){2}{\rule{0.400pt}{0.400pt}}
\multiput(187.00,121.59)(0.824,0.488){13}{\rule{0.750pt}{0.117pt}}
\multiput(187.00,120.17)(11.443,8.000){2}{\rule{0.375pt}{0.400pt}}
\multiput(200.00,129.59)(0.824,0.488){13}{\rule{0.750pt}{0.117pt}}
\multiput(200.00,128.17)(11.443,8.000){2}{\rule{0.375pt}{0.400pt}}
\multiput(213.00,137.59)(0.890,0.488){13}{\rule{0.800pt}{0.117pt}}
\multiput(213.00,136.17)(12.340,8.000){2}{\rule{0.400pt}{0.400pt}}
\multiput(227.00,145.59)(0.824,0.488){13}{\rule{0.750pt}{0.117pt}}
\multiput(227.00,144.17)(11.443,8.000){2}{\rule{0.375pt}{0.400pt}}
\multiput(240.00,153.59)(0.824,0.488){13}{\rule{0.750pt}{0.117pt}}
\multiput(240.00,152.17)(11.443,8.000){2}{\rule{0.375pt}{0.400pt}}
\multiput(253.00,161.59)(1.026,0.485){11}{\rule{0.900pt}{0.117pt}}
\multiput(253.00,160.17)(12.132,7.000){2}{\rule{0.450pt}{0.400pt}}
\multiput(267.00,168.59)(0.824,0.488){13}{\rule{0.750pt}{0.117pt}}
\multiput(267.00,167.17)(11.443,8.000){2}{\rule{0.375pt}{0.400pt}}
\multiput(280.00,176.59)(0.824,0.488){13}{\rule{0.750pt}{0.117pt}}
\multiput(280.00,175.17)(11.443,8.000){2}{\rule{0.375pt}{0.400pt}}
\multiput(293.00,184.59)(0.890,0.488){13}{\rule{0.800pt}{0.117pt}}
\multiput(293.00,183.17)(12.340,8.000){2}{\rule{0.400pt}{0.400pt}}
\multiput(307.00,192.59)(0.824,0.488){13}{\rule{0.750pt}{0.117pt}}
\multiput(307.00,191.17)(11.443,8.000){2}{\rule{0.375pt}{0.400pt}}
\multiput(320.00,200.59)(0.824,0.488){13}{\rule{0.750pt}{0.117pt}}
\multiput(320.00,199.17)(11.443,8.000){2}{\rule{0.375pt}{0.400pt}}
\multiput(333.00,208.59)(0.824,0.488){13}{\rule{0.750pt}{0.117pt}}
\multiput(333.00,207.17)(11.443,8.000){2}{\rule{0.375pt}{0.400pt}}
\multiput(346.00,216.59)(1.026,0.485){11}{\rule{0.900pt}{0.117pt}}
\multiput(346.00,215.17)(12.132,7.000){2}{\rule{0.450pt}{0.400pt}}
\multiput(360.00,223.59)(0.824,0.488){13}{\rule{0.750pt}{0.117pt}}
\multiput(360.00,222.17)(11.443,8.000){2}{\rule{0.375pt}{0.400pt}}
\multiput(373.00,231.59)(0.824,0.488){13}{\rule{0.750pt}{0.117pt}}
\multiput(373.00,230.17)(11.443,8.000){2}{\rule{0.375pt}{0.400pt}}
\multiput(386.00,239.59)(0.890,0.488){13}{\rule{0.800pt}{0.117pt}}
\multiput(386.00,238.17)(12.340,8.000){2}{\rule{0.400pt}{0.400pt}}
\multiput(400.00,247.59)(0.824,0.488){13}{\rule{0.750pt}{0.117pt}}
\multiput(400.00,246.17)(11.443,8.000){2}{\rule{0.375pt}{0.400pt}}
\multiput(413.00,255.59)(0.824,0.488){13}{\rule{0.750pt}{0.117pt}}
\multiput(413.00,254.17)(11.443,8.000){2}{\rule{0.375pt}{0.400pt}}
\multiput(426.00,263.59)(0.890,0.488){13}{\rule{0.800pt}{0.117pt}}
\multiput(426.00,262.17)(12.340,8.000){2}{\rule{0.400pt}{0.400pt}}
\multiput(440.00,271.59)(0.950,0.485){11}{\rule{0.843pt}{0.117pt}}
\multiput(440.00,270.17)(11.251,7.000){2}{\rule{0.421pt}{0.400pt}}
\multiput(453.00,278.59)(0.824,0.488){13}{\rule{0.750pt}{0.117pt}}
\multiput(453.00,277.17)(11.443,8.000){2}{\rule{0.375pt}{0.400pt}}
\multiput(466.00,286.59)(0.890,0.488){13}{\rule{0.800pt}{0.117pt}}
\multiput(466.00,285.17)(12.340,8.000){2}{\rule{0.400pt}{0.400pt}}
\multiput(480.00,294.59)(0.824,0.488){13}{\rule{0.750pt}{0.117pt}}
\multiput(480.00,293.17)(11.443,8.000){2}{\rule{0.375pt}{0.400pt}}
\multiput(493.00,302.59)(0.824,0.488){13}{\rule{0.750pt}{0.117pt}}
\multiput(493.00,301.17)(11.443,8.000){2}{\rule{0.375pt}{0.400pt}}
\multiput(506.00,310.59)(0.890,0.488){13}{\rule{0.800pt}{0.117pt}}
\multiput(506.00,309.17)(12.340,8.000){2}{\rule{0.400pt}{0.400pt}}
\multiput(520.00,318.59)(0.824,0.488){13}{\rule{0.750pt}{0.117pt}}
\multiput(520.00,317.17)(11.443,8.000){2}{\rule{0.375pt}{0.400pt}}
\multiput(533.00,326.59)(0.950,0.485){11}{\rule{0.843pt}{0.117pt}}
\multiput(533.00,325.17)(11.251,7.000){2}{\rule{0.421pt}{0.400pt}}
\multiput(546.00,333.59)(0.890,0.488){13}{\rule{0.800pt}{0.117pt}}
\multiput(546.00,332.17)(12.340,8.000){2}{\rule{0.400pt}{0.400pt}}
\multiput(560.00,341.59)(0.824,0.488){13}{\rule{0.750pt}{0.117pt}}
\multiput(560.00,340.17)(11.443,8.000){2}{\rule{0.375pt}{0.400pt}}
\multiput(573.00,349.59)(0.824,0.488){13}{\rule{0.750pt}{0.117pt}}
\multiput(573.00,348.17)(11.443,8.000){2}{\rule{0.375pt}{0.400pt}}
\multiput(586.00,357.59)(0.890,0.488){13}{\rule{0.800pt}{0.117pt}}
\multiput(586.00,356.17)(12.340,8.000){2}{\rule{0.400pt}{0.400pt}}
\multiput(600.00,365.59)(0.824,0.488){13}{\rule{0.750pt}{0.117pt}}
\multiput(600.00,364.17)(11.443,8.000){2}{\rule{0.375pt}{0.400pt}}
\multiput(613.00,373.59)(0.824,0.488){13}{\rule{0.750pt}{0.117pt}}
\multiput(613.00,372.17)(11.443,8.000){2}{\rule{0.375pt}{0.400pt}}
\multiput(626.00,381.59)(1.026,0.485){11}{\rule{0.900pt}{0.117pt}}
\multiput(626.00,380.17)(12.132,7.000){2}{\rule{0.450pt}{0.400pt}}
\multiput(640.00,388.59)(0.824,0.488){13}{\rule{0.750pt}{0.117pt}}
\multiput(640.00,387.17)(11.443,8.000){2}{\rule{0.375pt}{0.400pt}}
\multiput(653.00,396.59)(0.824,0.488){13}{\rule{0.750pt}{0.117pt}}
\multiput(653.00,395.17)(11.443,8.000){2}{\rule{0.375pt}{0.400pt}}
\multiput(666.00,404.59)(0.890,0.488){13}{\rule{0.800pt}{0.117pt}}
\multiput(666.00,403.17)(12.340,8.000){2}{\rule{0.400pt}{0.400pt}}
\multiput(680.00,412.59)(0.824,0.488){13}{\rule{0.750pt}{0.117pt}}
\multiput(680.00,411.17)(11.443,8.000){2}{\rule{0.375pt}{0.400pt}}
\multiput(693.00,420.59)(0.824,0.488){13}{\rule{0.750pt}{0.117pt}}
\multiput(693.00,419.17)(11.443,8.000){2}{\rule{0.375pt}{0.400pt}}
\multiput(706.00,428.59)(0.890,0.488){13}{\rule{0.800pt}{0.117pt}}
\multiput(706.00,427.17)(12.340,8.000){2}{\rule{0.400pt}{0.400pt}}
\multiput(720.00,436.59)(0.950,0.485){11}{\rule{0.843pt}{0.117pt}}
\multiput(720.00,435.17)(11.251,7.000){2}{\rule{0.421pt}{0.400pt}}
\multiput(733.00,443.59)(0.824,0.488){13}{\rule{0.750pt}{0.117pt}}
\multiput(733.00,442.17)(11.443,8.000){2}{\rule{0.375pt}{0.400pt}}
\multiput(746.00,451.59)(0.890,0.488){13}{\rule{0.800pt}{0.117pt}}
\multiput(746.00,450.17)(12.340,8.000){2}{\rule{0.400pt}{0.400pt}}
\multiput(760.00,459.59)(0.824,0.488){13}{\rule{0.750pt}{0.117pt}}
\multiput(760.00,458.17)(11.443,8.000){2}{\rule{0.375pt}{0.400pt}}
\multiput(773.00,467.59)(0.824,0.488){13}{\rule{0.750pt}{0.117pt}}
\multiput(773.00,466.17)(11.443,8.000){2}{\rule{0.375pt}{0.400pt}}
\multiput(786.00,475.59)(0.824,0.488){13}{\rule{0.750pt}{0.117pt}}
\multiput(786.00,474.17)(11.443,8.000){2}{\rule{0.375pt}{0.400pt}}
\multiput(799.00,483.59)(0.890,0.488){13}{\rule{0.800pt}{0.117pt}}
\multiput(799.00,482.17)(12.340,8.000){2}{\rule{0.400pt}{0.400pt}}
\multiput(813.00,491.59)(0.824,0.488){13}{\rule{0.750pt}{0.117pt}}
\multiput(813.00,490.17)(11.443,8.000){2}{\rule{0.375pt}{0.400pt}}
\multiput(826.00,499.59)(0.950,0.485){11}{\rule{0.843pt}{0.117pt}}
\multiput(826.00,498.17)(11.251,7.000){2}{\rule{0.421pt}{0.400pt}}
\multiput(839.00,506.59)(0.890,0.488){13}{\rule{0.800pt}{0.117pt}}
\multiput(839.00,505.17)(12.340,8.000){2}{\rule{0.400pt}{0.400pt}}
\multiput(853.00,514.59)(0.824,0.488){13}{\rule{0.750pt}{0.117pt}}
\multiput(853.00,513.17)(11.443,8.000){2}{\rule{0.375pt}{0.400pt}}
\multiput(866.00,522.59)(0.824,0.488){13}{\rule{0.750pt}{0.117pt}}
\multiput(866.00,521.17)(11.443,8.000){2}{\rule{0.375pt}{0.400pt}}
\multiput(879.00,530.59)(0.890,0.488){13}{\rule{0.800pt}{0.117pt}}
\multiput(879.00,529.17)(12.340,8.000){2}{\rule{0.400pt}{0.400pt}}
\multiput(893.00,538.59)(0.824,0.488){13}{\rule{0.750pt}{0.117pt}}
\multiput(893.00,537.17)(11.443,8.000){2}{\rule{0.375pt}{0.400pt}}
\multiput(906.00,546.59)(0.824,0.488){13}{\rule{0.750pt}{0.117pt}}
\multiput(906.00,545.17)(11.443,8.000){2}{\rule{0.375pt}{0.400pt}}
\multiput(919.00,554.59)(1.026,0.485){11}{\rule{0.900pt}{0.117pt}}
\multiput(919.00,553.17)(12.132,7.000){2}{\rule{0.450pt}{0.400pt}}
\multiput(933.00,561.59)(0.824,0.488){13}{\rule{0.750pt}{0.117pt}}
\multiput(933.00,560.17)(11.443,8.000){2}{\rule{0.375pt}{0.400pt}}
\multiput(946.00,569.59)(0.824,0.488){13}{\rule{0.750pt}{0.117pt}}
\multiput(946.00,568.17)(11.443,8.000){2}{\rule{0.375pt}{0.400pt}}
\multiput(959.00,577.59)(0.890,0.488){13}{\rule{0.800pt}{0.117pt}}
\multiput(959.00,576.17)(12.340,8.000){2}{\rule{0.400pt}{0.400pt}}
\multiput(973.00,585.59)(0.824,0.488){13}{\rule{0.750pt}{0.117pt}}
\multiput(973.00,584.17)(11.443,8.000){2}{\rule{0.375pt}{0.400pt}}
\multiput(986.00,593.59)(0.824,0.488){13}{\rule{0.750pt}{0.117pt}}
\multiput(986.00,592.17)(11.443,8.000){2}{\rule{0.375pt}{0.400pt}}
\multiput(999.00,601.59)(0.890,0.488){13}{\rule{0.800pt}{0.117pt}}
\multiput(999.00,600.17)(12.340,8.000){2}{\rule{0.400pt}{0.400pt}}
\multiput(1013.00,609.59)(0.950,0.485){11}{\rule{0.843pt}{0.117pt}}
\multiput(1013.00,608.17)(11.251,7.000){2}{\rule{0.421pt}{0.400pt}}
\multiput(1026.00,616.59)(0.824,0.488){13}{\rule{0.750pt}{0.117pt}}
\multiput(1026.00,615.17)(11.443,8.000){2}{\rule{0.375pt}{0.400pt}}
\multiput(1039.00,624.59)(0.890,0.488){13}{\rule{0.800pt}{0.117pt}}
\multiput(1039.00,623.17)(12.340,8.000){2}{\rule{0.400pt}{0.400pt}}
\multiput(1053.00,632.59)(0.824,0.488){13}{\rule{0.750pt}{0.117pt}}
\multiput(1053.00,631.17)(11.443,8.000){2}{\rule{0.375pt}{0.400pt}}
\multiput(1066.00,640.59)(0.824,0.488){13}{\rule{0.750pt}{0.117pt}}
\multiput(1066.00,639.17)(11.443,8.000){2}{\rule{0.375pt}{0.400pt}}
\multiput(1079.00,648.59)(0.890,0.488){13}{\rule{0.800pt}{0.117pt}}
\multiput(1079.00,647.17)(12.340,8.000){2}{\rule{0.400pt}{0.400pt}}
\multiput(1093.00,656.59)(0.824,0.488){13}{\rule{0.750pt}{0.117pt}}
\multiput(1093.00,655.17)(11.443,8.000){2}{\rule{0.375pt}{0.400pt}}
\multiput(1106.00,664.59)(0.950,0.485){11}{\rule{0.843pt}{0.117pt}}
\multiput(1106.00,663.17)(11.251,7.000){2}{\rule{0.421pt}{0.400pt}}
\multiput(1119.00,671.59)(0.890,0.488){13}{\rule{0.800pt}{0.117pt}}
\multiput(1119.00,670.17)(12.340,8.000){2}{\rule{0.400pt}{0.400pt}}
\multiput(1133.00,679.59)(0.824,0.488){13}{\rule{0.750pt}{0.117pt}}
\multiput(1133.00,678.17)(11.443,8.000){2}{\rule{0.375pt}{0.400pt}}
\multiput(1146.00,687.59)(0.824,0.488){13}{\rule{0.750pt}{0.117pt}}
\multiput(1146.00,686.17)(11.443,8.000){2}{\rule{0.375pt}{0.400pt}}
\multiput(1159.00,695.59)(0.890,0.488){13}{\rule{0.800pt}{0.117pt}}
\multiput(1159.00,694.17)(12.340,8.000){2}{\rule{0.400pt}{0.400pt}}
\multiput(1173.00,703.59)(0.824,0.488){13}{\rule{0.750pt}{0.117pt}}
\multiput(1173.00,702.17)(11.443,8.000){2}{\rule{0.375pt}{0.400pt}}
\multiput(1186.00,711.59)(0.824,0.488){13}{\rule{0.750pt}{0.117pt}}
\multiput(1186.00,710.17)(11.443,8.000){2}{\rule{0.375pt}{0.400pt}}
\multiput(1199.00,719.59)(1.026,0.485){11}{\rule{0.900pt}{0.117pt}}
\multiput(1199.00,718.17)(12.132,7.000){2}{\rule{0.450pt}{0.400pt}}
\multiput(1213.00,726.59)(0.824,0.488){13}{\rule{0.750pt}{0.117pt}}
\multiput(1213.00,725.17)(11.443,8.000){2}{\rule{0.375pt}{0.400pt}}
\multiput(1226.00,734.59)(0.824,0.488){13}{\rule{0.750pt}{0.117pt}}
\multiput(1226.00,733.17)(11.443,8.000){2}{\rule{0.375pt}{0.400pt}}
\multiput(1239.00,742.59)(0.824,0.488){13}{\rule{0.750pt}{0.117pt}}
\multiput(1239.00,741.17)(11.443,8.000){2}{\rule{0.375pt}{0.400pt}}
\multiput(1252.00,750.59)(0.890,0.488){13}{\rule{0.800pt}{0.117pt}}
\multiput(1252.00,749.17)(12.340,8.000){2}{\rule{0.400pt}{0.400pt}}
\multiput(1266.00,758.59)(0.824,0.488){13}{\rule{0.750pt}{0.117pt}}
\multiput(1266.00,757.17)(11.443,8.000){2}{\rule{0.375pt}{0.400pt}}
\multiput(1279.00,766.59)(0.824,0.488){13}{\rule{0.750pt}{0.117pt}}
\multiput(1279.00,765.17)(11.443,8.000){2}{\rule{0.375pt}{0.400pt}}
\multiput(1292.00,774.59)(1.026,0.485){11}{\rule{0.900pt}{0.117pt}}
\multiput(1292.00,773.17)(12.132,7.000){2}{\rule{0.450pt}{0.400pt}}
\multiput(1306.00,781.59)(0.824,0.488){13}{\rule{0.750pt}{0.117pt}}
\multiput(1306.00,780.17)(11.443,8.000){2}{\rule{0.375pt}{0.400pt}}
\multiput(1319.00,789.59)(0.824,0.488){13}{\rule{0.750pt}{0.117pt}}
\multiput(1319.00,788.17)(11.443,8.000){2}{\rule{0.375pt}{0.400pt}}
\multiput(1332.00,797.59)(0.890,0.488){13}{\rule{0.800pt}{0.117pt}}
\multiput(1332.00,796.17)(12.340,8.000){2}{\rule{0.400pt}{0.400pt}}
\multiput(1346.00,805.59)(0.824,0.488){13}{\rule{0.750pt}{0.117pt}}
\multiput(1346.00,804.17)(11.443,8.000){2}{\rule{0.375pt}{0.400pt}}
\multiput(1359.00,813.59)(0.824,0.488){13}{\rule{0.750pt}{0.117pt}}
\multiput(1359.00,812.17)(11.443,8.000){2}{\rule{0.375pt}{0.400pt}}
\multiput(1372.00,821.59)(0.890,0.488){13}{\rule{0.800pt}{0.117pt}}
\multiput(1372.00,820.17)(12.340,8.000){2}{\rule{0.400pt}{0.400pt}}
\multiput(1386.00,829.59)(0.950,0.485){11}{\rule{0.843pt}{0.117pt}}
\multiput(1386.00,828.17)(11.251,7.000){2}{\rule{0.421pt}{0.400pt}}
\multiput(1399.00,836.59)(0.824,0.488){13}{\rule{0.750pt}{0.117pt}}
\multiput(1399.00,835.17)(11.443,8.000){2}{\rule{0.375pt}{0.400pt}}
\multiput(1412.00,844.59)(0.890,0.488){13}{\rule{0.800pt}{0.117pt}}
\multiput(1412.00,843.17)(12.340,8.000){2}{\rule{0.400pt}{0.400pt}}
\multiput(1426.00,852.59)(0.824,0.488){13}{\rule{0.750pt}{0.117pt}}
\multiput(1426.00,851.17)(11.443,8.000){2}{\rule{0.375pt}{0.400pt}}
\put(120,82){\usebox{\plotpoint}}
\multiput(120.58,82.00)(0.493,0.616){23}{\rule{0.119pt}{0.592pt}}
\multiput(119.17,82.00)(13.000,14.771){2}{\rule{0.400pt}{0.296pt}}
\multiput(133.58,98.00)(0.494,0.534){25}{\rule{0.119pt}{0.529pt}}
\multiput(132.17,98.00)(14.000,13.903){2}{\rule{0.400pt}{0.264pt}}
\multiput(147.58,113.00)(0.493,0.576){23}{\rule{0.119pt}{0.562pt}}
\multiput(146.17,113.00)(13.000,13.834){2}{\rule{0.400pt}{0.281pt}}
\multiput(160.58,128.00)(0.493,0.616){23}{\rule{0.119pt}{0.592pt}}
\multiput(159.17,128.00)(13.000,14.771){2}{\rule{0.400pt}{0.296pt}}
\multiput(173.58,144.00)(0.494,0.534){25}{\rule{0.119pt}{0.529pt}}
\multiput(172.17,144.00)(14.000,13.903){2}{\rule{0.400pt}{0.264pt}}
\multiput(187.58,159.00)(0.493,0.536){23}{\rule{0.119pt}{0.531pt}}
\multiput(186.17,159.00)(13.000,12.898){2}{\rule{0.400pt}{0.265pt}}
\multiput(200.58,173.00)(0.493,0.576){23}{\rule{0.119pt}{0.562pt}}
\multiput(199.17,173.00)(13.000,13.834){2}{\rule{0.400pt}{0.281pt}}
\multiput(213.58,188.00)(0.494,0.534){25}{\rule{0.119pt}{0.529pt}}
\multiput(212.17,188.00)(14.000,13.903){2}{\rule{0.400pt}{0.264pt}}
\multiput(227.58,203.00)(0.493,0.536){23}{\rule{0.119pt}{0.531pt}}
\multiput(226.17,203.00)(13.000,12.898){2}{\rule{0.400pt}{0.265pt}}
\multiput(240.58,217.00)(0.493,0.536){23}{\rule{0.119pt}{0.531pt}}
\multiput(239.17,217.00)(13.000,12.898){2}{\rule{0.400pt}{0.265pt}}
\multiput(253.00,231.58)(0.497,0.494){25}{\rule{0.500pt}{0.119pt}}
\multiput(253.00,230.17)(12.962,14.000){2}{\rule{0.250pt}{0.400pt}}
\multiput(267.58,245.00)(0.493,0.536){23}{\rule{0.119pt}{0.531pt}}
\multiput(266.17,245.00)(13.000,12.898){2}{\rule{0.400pt}{0.265pt}}
\multiput(280.58,259.00)(0.493,0.536){23}{\rule{0.119pt}{0.531pt}}
\multiput(279.17,259.00)(13.000,12.898){2}{\rule{0.400pt}{0.265pt}}
\multiput(293.00,273.58)(0.536,0.493){23}{\rule{0.531pt}{0.119pt}}
\multiput(293.00,272.17)(12.898,13.000){2}{\rule{0.265pt}{0.400pt}}
\multiput(307.58,286.00)(0.493,0.536){23}{\rule{0.119pt}{0.531pt}}
\multiput(306.17,286.00)(13.000,12.898){2}{\rule{0.400pt}{0.265pt}}
\multiput(320.00,300.58)(0.497,0.493){23}{\rule{0.500pt}{0.119pt}}
\multiput(320.00,299.17)(11.962,13.000){2}{\rule{0.250pt}{0.400pt}}
\multiput(333.00,313.58)(0.497,0.493){23}{\rule{0.500pt}{0.119pt}}
\multiput(333.00,312.17)(11.962,13.000){2}{\rule{0.250pt}{0.400pt}}
\multiput(346.00,326.58)(0.536,0.493){23}{\rule{0.531pt}{0.119pt}}
\multiput(346.00,325.17)(12.898,13.000){2}{\rule{0.265pt}{0.400pt}}
\multiput(360.00,339.58)(0.497,0.493){23}{\rule{0.500pt}{0.119pt}}
\multiput(360.00,338.17)(11.962,13.000){2}{\rule{0.250pt}{0.400pt}}
\multiput(373.00,352.58)(0.497,0.493){23}{\rule{0.500pt}{0.119pt}}
\multiput(373.00,351.17)(11.962,13.000){2}{\rule{0.250pt}{0.400pt}}
\multiput(386.00,365.58)(0.582,0.492){21}{\rule{0.567pt}{0.119pt}}
\multiput(386.00,364.17)(12.824,12.000){2}{\rule{0.283pt}{0.400pt}}
\multiput(400.00,377.58)(0.539,0.492){21}{\rule{0.533pt}{0.119pt}}
\multiput(400.00,376.17)(11.893,12.000){2}{\rule{0.267pt}{0.400pt}}
\multiput(413.00,389.58)(0.497,0.493){23}{\rule{0.500pt}{0.119pt}}
\multiput(413.00,388.17)(11.962,13.000){2}{\rule{0.250pt}{0.400pt}}
\multiput(426.00,402.58)(0.637,0.492){19}{\rule{0.609pt}{0.118pt}}
\multiput(426.00,401.17)(12.736,11.000){2}{\rule{0.305pt}{0.400pt}}
\multiput(440.00,413.58)(0.539,0.492){21}{\rule{0.533pt}{0.119pt}}
\multiput(440.00,412.17)(11.893,12.000){2}{\rule{0.267pt}{0.400pt}}
\multiput(453.00,425.58)(0.539,0.492){21}{\rule{0.533pt}{0.119pt}}
\multiput(453.00,424.17)(11.893,12.000){2}{\rule{0.267pt}{0.400pt}}
\multiput(466.00,437.58)(0.637,0.492){19}{\rule{0.609pt}{0.118pt}}
\multiput(466.00,436.17)(12.736,11.000){2}{\rule{0.305pt}{0.400pt}}
\multiput(480.00,448.58)(0.539,0.492){21}{\rule{0.533pt}{0.119pt}}
\multiput(480.00,447.17)(11.893,12.000){2}{\rule{0.267pt}{0.400pt}}
\multiput(493.00,460.58)(0.590,0.492){19}{\rule{0.573pt}{0.118pt}}
\multiput(493.00,459.17)(11.811,11.000){2}{\rule{0.286pt}{0.400pt}}
\multiput(506.00,471.58)(0.637,0.492){19}{\rule{0.609pt}{0.118pt}}
\multiput(506.00,470.17)(12.736,11.000){2}{\rule{0.305pt}{0.400pt}}
\multiput(520.00,482.58)(0.590,0.492){19}{\rule{0.573pt}{0.118pt}}
\multiput(520.00,481.17)(11.811,11.000){2}{\rule{0.286pt}{0.400pt}}
\multiput(533.00,493.58)(0.590,0.492){19}{\rule{0.573pt}{0.118pt}}
\multiput(533.00,492.17)(11.811,11.000){2}{\rule{0.286pt}{0.400pt}}
\multiput(546.00,504.58)(0.704,0.491){17}{\rule{0.660pt}{0.118pt}}
\multiput(546.00,503.17)(12.630,10.000){2}{\rule{0.330pt}{0.400pt}}
\multiput(560.00,514.58)(0.590,0.492){19}{\rule{0.573pt}{0.118pt}}
\multiput(560.00,513.17)(11.811,11.000){2}{\rule{0.286pt}{0.400pt}}
\multiput(573.00,525.58)(0.652,0.491){17}{\rule{0.620pt}{0.118pt}}
\multiput(573.00,524.17)(11.713,10.000){2}{\rule{0.310pt}{0.400pt}}
\multiput(586.00,535.58)(0.704,0.491){17}{\rule{0.660pt}{0.118pt}}
\multiput(586.00,534.17)(12.630,10.000){2}{\rule{0.330pt}{0.400pt}}
\multiput(600.00,545.58)(0.652,0.491){17}{\rule{0.620pt}{0.118pt}}
\multiput(600.00,544.17)(11.713,10.000){2}{\rule{0.310pt}{0.400pt}}
\multiput(613.00,555.58)(0.652,0.491){17}{\rule{0.620pt}{0.118pt}}
\multiput(613.00,554.17)(11.713,10.000){2}{\rule{0.310pt}{0.400pt}}
\multiput(626.00,565.59)(0.786,0.489){15}{\rule{0.722pt}{0.118pt}}
\multiput(626.00,564.17)(12.501,9.000){2}{\rule{0.361pt}{0.400pt}}
\multiput(640.00,574.58)(0.652,0.491){17}{\rule{0.620pt}{0.118pt}}
\multiput(640.00,573.17)(11.713,10.000){2}{\rule{0.310pt}{0.400pt}}
\multiput(653.00,584.59)(0.728,0.489){15}{\rule{0.678pt}{0.118pt}}
\multiput(653.00,583.17)(11.593,9.000){2}{\rule{0.339pt}{0.400pt}}
\multiput(666.00,593.59)(0.786,0.489){15}{\rule{0.722pt}{0.118pt}}
\multiput(666.00,592.17)(12.501,9.000){2}{\rule{0.361pt}{0.400pt}}
\multiput(680.00,602.59)(0.728,0.489){15}{\rule{0.678pt}{0.118pt}}
\multiput(680.00,601.17)(11.593,9.000){2}{\rule{0.339pt}{0.400pt}}
\multiput(693.00,611.59)(0.728,0.489){15}{\rule{0.678pt}{0.118pt}}
\multiput(693.00,610.17)(11.593,9.000){2}{\rule{0.339pt}{0.400pt}}
\multiput(706.00,620.59)(0.786,0.489){15}{\rule{0.722pt}{0.118pt}}
\multiput(706.00,619.17)(12.501,9.000){2}{\rule{0.361pt}{0.400pt}}
\multiput(720.00,629.59)(0.824,0.488){13}{\rule{0.750pt}{0.117pt}}
\multiput(720.00,628.17)(11.443,8.000){2}{\rule{0.375pt}{0.400pt}}
\multiput(733.00,637.59)(0.824,0.488){13}{\rule{0.750pt}{0.117pt}}
\multiput(733.00,636.17)(11.443,8.000){2}{\rule{0.375pt}{0.400pt}}
\multiput(746.00,645.59)(0.786,0.489){15}{\rule{0.722pt}{0.118pt}}
\multiput(746.00,644.17)(12.501,9.000){2}{\rule{0.361pt}{0.400pt}}
\multiput(760.00,654.59)(0.824,0.488){13}{\rule{0.750pt}{0.117pt}}
\multiput(760.00,653.17)(11.443,8.000){2}{\rule{0.375pt}{0.400pt}}
\multiput(773.00,662.59)(0.950,0.485){11}{\rule{0.843pt}{0.117pt}}
\multiput(773.00,661.17)(11.251,7.000){2}{\rule{0.421pt}{0.400pt}}
\multiput(786.00,669.59)(0.824,0.488){13}{\rule{0.750pt}{0.117pt}}
\multiput(786.00,668.17)(11.443,8.000){2}{\rule{0.375pt}{0.400pt}}
\multiput(799.00,677.59)(0.890,0.488){13}{\rule{0.800pt}{0.117pt}}
\multiput(799.00,676.17)(12.340,8.000){2}{\rule{0.400pt}{0.400pt}}
\multiput(813.00,685.59)(0.950,0.485){11}{\rule{0.843pt}{0.117pt}}
\multiput(813.00,684.17)(11.251,7.000){2}{\rule{0.421pt}{0.400pt}}
\multiput(826.00,692.59)(0.950,0.485){11}{\rule{0.843pt}{0.117pt}}
\multiput(826.00,691.17)(11.251,7.000){2}{\rule{0.421pt}{0.400pt}}
\multiput(839.00,699.59)(1.026,0.485){11}{\rule{0.900pt}{0.117pt}}
\multiput(839.00,698.17)(12.132,7.000){2}{\rule{0.450pt}{0.400pt}}
\multiput(853.00,706.59)(0.950,0.485){11}{\rule{0.843pt}{0.117pt}}
\multiput(853.00,705.17)(11.251,7.000){2}{\rule{0.421pt}{0.400pt}}
\multiput(866.00,713.59)(0.950,0.485){11}{\rule{0.843pt}{0.117pt}}
\multiput(866.00,712.17)(11.251,7.000){2}{\rule{0.421pt}{0.400pt}}
\multiput(879.00,720.59)(1.026,0.485){11}{\rule{0.900pt}{0.117pt}}
\multiput(879.00,719.17)(12.132,7.000){2}{\rule{0.450pt}{0.400pt}}
\multiput(893.00,727.59)(1.123,0.482){9}{\rule{0.967pt}{0.116pt}}
\multiput(893.00,726.17)(10.994,6.000){2}{\rule{0.483pt}{0.400pt}}
\multiput(906.00,733.59)(1.123,0.482){9}{\rule{0.967pt}{0.116pt}}
\multiput(906.00,732.17)(10.994,6.000){2}{\rule{0.483pt}{0.400pt}}
\multiput(919.00,739.59)(1.214,0.482){9}{\rule{1.033pt}{0.116pt}}
\multiput(919.00,738.17)(11.855,6.000){2}{\rule{0.517pt}{0.400pt}}
\multiput(933.00,745.59)(1.123,0.482){9}{\rule{0.967pt}{0.116pt}}
\multiput(933.00,744.17)(10.994,6.000){2}{\rule{0.483pt}{0.400pt}}
\multiput(946.00,751.59)(1.123,0.482){9}{\rule{0.967pt}{0.116pt}}
\multiput(946.00,750.17)(10.994,6.000){2}{\rule{0.483pt}{0.400pt}}
\multiput(959.00,757.59)(1.214,0.482){9}{\rule{1.033pt}{0.116pt}}
\multiput(959.00,756.17)(11.855,6.000){2}{\rule{0.517pt}{0.400pt}}
\multiput(973.00,763.59)(1.378,0.477){7}{\rule{1.140pt}{0.115pt}}
\multiput(973.00,762.17)(10.634,5.000){2}{\rule{0.570pt}{0.400pt}}
\multiput(986.00,768.59)(1.123,0.482){9}{\rule{0.967pt}{0.116pt}}
\multiput(986.00,767.17)(10.994,6.000){2}{\rule{0.483pt}{0.400pt}}
\multiput(999.00,774.59)(1.489,0.477){7}{\rule{1.220pt}{0.115pt}}
\multiput(999.00,773.17)(11.468,5.000){2}{\rule{0.610pt}{0.400pt}}
\multiput(1013.00,779.59)(1.378,0.477){7}{\rule{1.140pt}{0.115pt}}
\multiput(1013.00,778.17)(10.634,5.000){2}{\rule{0.570pt}{0.400pt}}
\multiput(1026.00,784.59)(1.378,0.477){7}{\rule{1.140pt}{0.115pt}}
\multiput(1026.00,783.17)(10.634,5.000){2}{\rule{0.570pt}{0.400pt}}
\multiput(1039.00,789.60)(1.943,0.468){5}{\rule{1.500pt}{0.113pt}}
\multiput(1039.00,788.17)(10.887,4.000){2}{\rule{0.750pt}{0.400pt}}
\multiput(1053.00,793.59)(1.378,0.477){7}{\rule{1.140pt}{0.115pt}}
\multiput(1053.00,792.17)(10.634,5.000){2}{\rule{0.570pt}{0.400pt}}
\multiput(1066.00,798.60)(1.797,0.468){5}{\rule{1.400pt}{0.113pt}}
\multiput(1066.00,797.17)(10.094,4.000){2}{\rule{0.700pt}{0.400pt}}
\multiput(1079.00,802.60)(1.943,0.468){5}{\rule{1.500pt}{0.113pt}}
\multiput(1079.00,801.17)(10.887,4.000){2}{\rule{0.750pt}{0.400pt}}
\multiput(1093.00,806.60)(1.797,0.468){5}{\rule{1.400pt}{0.113pt}}
\multiput(1093.00,805.17)(10.094,4.000){2}{\rule{0.700pt}{0.400pt}}
\multiput(1106.00,810.60)(1.797,0.468){5}{\rule{1.400pt}{0.113pt}}
\multiput(1106.00,809.17)(10.094,4.000){2}{\rule{0.700pt}{0.400pt}}
\multiput(1119.00,814.60)(1.943,0.468){5}{\rule{1.500pt}{0.113pt}}
\multiput(1119.00,813.17)(10.887,4.000){2}{\rule{0.750pt}{0.400pt}}
\multiput(1133.00,818.60)(1.797,0.468){5}{\rule{1.400pt}{0.113pt}}
\multiput(1133.00,817.17)(10.094,4.000){2}{\rule{0.700pt}{0.400pt}}
\multiput(1146.00,822.61)(2.695,0.447){3}{\rule{1.833pt}{0.108pt}}
\multiput(1146.00,821.17)(9.195,3.000){2}{\rule{0.917pt}{0.400pt}}
\multiput(1159.00,825.61)(2.918,0.447){3}{\rule{1.967pt}{0.108pt}}
\multiput(1159.00,824.17)(9.918,3.000){2}{\rule{0.983pt}{0.400pt}}
\multiput(1173.00,828.61)(2.695,0.447){3}{\rule{1.833pt}{0.108pt}}
\multiput(1173.00,827.17)(9.195,3.000){2}{\rule{0.917pt}{0.400pt}}
\multiput(1186.00,831.61)(2.695,0.447){3}{\rule{1.833pt}{0.108pt}}
\multiput(1186.00,830.17)(9.195,3.000){2}{\rule{0.917pt}{0.400pt}}
\multiput(1199.00,834.61)(2.918,0.447){3}{\rule{1.967pt}{0.108pt}}
\multiput(1199.00,833.17)(9.918,3.000){2}{\rule{0.983pt}{0.400pt}}
\multiput(1213.00,837.61)(2.695,0.447){3}{\rule{1.833pt}{0.108pt}}
\multiput(1213.00,836.17)(9.195,3.000){2}{\rule{0.917pt}{0.400pt}}
\put(1226,840.17){\rule{2.700pt}{0.400pt}}
\multiput(1226.00,839.17)(7.396,2.000){2}{\rule{1.350pt}{0.400pt}}
\put(1239,842.17){\rule{2.700pt}{0.400pt}}
\multiput(1239.00,841.17)(7.396,2.000){2}{\rule{1.350pt}{0.400pt}}
\multiput(1252.00,844.61)(2.918,0.447){3}{\rule{1.967pt}{0.108pt}}
\multiput(1252.00,843.17)(9.918,3.000){2}{\rule{0.983pt}{0.400pt}}
\put(1266,847.17){\rule{2.700pt}{0.400pt}}
\multiput(1266.00,846.17)(7.396,2.000){2}{\rule{1.350pt}{0.400pt}}
\put(1279,848.67){\rule{3.132pt}{0.400pt}}
\multiput(1279.00,848.17)(6.500,1.000){2}{\rule{1.566pt}{0.400pt}}
\put(1292,850.17){\rule{2.900pt}{0.400pt}}
\multiput(1292.00,849.17)(7.981,2.000){2}{\rule{1.450pt}{0.400pt}}
\put(1306,852.17){\rule{2.700pt}{0.400pt}}
\multiput(1306.00,851.17)(7.396,2.000){2}{\rule{1.350pt}{0.400pt}}
\put(1319,853.67){\rule{3.132pt}{0.400pt}}
\multiput(1319.00,853.17)(6.500,1.000){2}{\rule{1.566pt}{0.400pt}}
\put(1332,854.67){\rule{3.373pt}{0.400pt}}
\multiput(1332.00,854.17)(7.000,1.000){2}{\rule{1.686pt}{0.400pt}}
\put(1346,855.67){\rule{3.132pt}{0.400pt}}
\multiput(1346.00,855.17)(6.500,1.000){2}{\rule{1.566pt}{0.400pt}}
\put(1359,856.67){\rule{3.132pt}{0.400pt}}
\multiput(1359.00,856.17)(6.500,1.000){2}{\rule{1.566pt}{0.400pt}}
\put(1372,857.67){\rule{3.373pt}{0.400pt}}
\multiput(1372.00,857.17)(7.000,1.000){2}{\rule{1.686pt}{0.400pt}}
\put(1399,858.67){\rule{3.132pt}{0.400pt}}
\multiput(1399.00,858.17)(6.500,1.000){2}{\rule{1.566pt}{0.400pt}}
\put(1386.0,859.0){\rule[-0.200pt]{3.132pt}{0.400pt}}
\put(1412.0,860.0){\rule[-0.200pt]{6.504pt}{0.400pt}}
\put(120,82){\usebox{\plotpoint}}
\multiput(120.58,82.00)(0.493,0.616){23}{\rule{0.119pt}{0.592pt}}
\multiput(119.17,82.00)(13.000,14.771){2}{\rule{0.400pt}{0.296pt}}
\multiput(133.58,98.00)(0.494,0.534){25}{\rule{0.119pt}{0.529pt}}
\multiput(132.17,98.00)(14.000,13.903){2}{\rule{0.400pt}{0.264pt}}
\multiput(147.58,113.00)(0.493,0.576){23}{\rule{0.119pt}{0.562pt}}
\multiput(146.17,113.00)(13.000,13.834){2}{\rule{0.400pt}{0.281pt}}
\multiput(160.58,128.00)(0.493,0.536){23}{\rule{0.119pt}{0.531pt}}
\multiput(159.17,128.00)(13.000,12.898){2}{\rule{0.400pt}{0.265pt}}
\multiput(173.58,142.00)(0.494,0.534){25}{\rule{0.119pt}{0.529pt}}
\multiput(172.17,142.00)(14.000,13.903){2}{\rule{0.400pt}{0.264pt}}
\multiput(187.58,157.00)(0.493,0.536){23}{\rule{0.119pt}{0.531pt}}
\multiput(186.17,157.00)(13.000,12.898){2}{\rule{0.400pt}{0.265pt}}
\multiput(200.00,171.58)(0.497,0.493){23}{\rule{0.500pt}{0.119pt}}
\multiput(200.00,170.17)(11.962,13.000){2}{\rule{0.250pt}{0.400pt}}
\multiput(213.00,184.58)(0.497,0.494){25}{\rule{0.500pt}{0.119pt}}
\multiput(213.00,183.17)(12.962,14.000){2}{\rule{0.250pt}{0.400pt}}
\multiput(227.00,198.58)(0.497,0.493){23}{\rule{0.500pt}{0.119pt}}
\multiput(227.00,197.17)(11.962,13.000){2}{\rule{0.250pt}{0.400pt}}
\multiput(240.00,211.58)(0.539,0.492){21}{\rule{0.533pt}{0.119pt}}
\multiput(240.00,210.17)(11.893,12.000){2}{\rule{0.267pt}{0.400pt}}
\multiput(253.00,223.58)(0.536,0.493){23}{\rule{0.531pt}{0.119pt}}
\multiput(253.00,222.17)(12.898,13.000){2}{\rule{0.265pt}{0.400pt}}
\multiput(267.00,236.58)(0.539,0.492){21}{\rule{0.533pt}{0.119pt}}
\multiput(267.00,235.17)(11.893,12.000){2}{\rule{0.267pt}{0.400pt}}
\multiput(280.00,248.58)(0.590,0.492){19}{\rule{0.573pt}{0.118pt}}
\multiput(280.00,247.17)(11.811,11.000){2}{\rule{0.286pt}{0.400pt}}
\multiput(293.00,259.58)(0.582,0.492){21}{\rule{0.567pt}{0.119pt}}
\multiput(293.00,258.17)(12.824,12.000){2}{\rule{0.283pt}{0.400pt}}
\multiput(307.00,271.58)(0.590,0.492){19}{\rule{0.573pt}{0.118pt}}
\multiput(307.00,270.17)(11.811,11.000){2}{\rule{0.286pt}{0.400pt}}
\multiput(320.00,282.58)(0.590,0.492){19}{\rule{0.573pt}{0.118pt}}
\multiput(320.00,281.17)(11.811,11.000){2}{\rule{0.286pt}{0.400pt}}
\multiput(333.00,293.58)(0.652,0.491){17}{\rule{0.620pt}{0.118pt}}
\multiput(333.00,292.17)(11.713,10.000){2}{\rule{0.310pt}{0.400pt}}
\multiput(346.00,303.58)(0.704,0.491){17}{\rule{0.660pt}{0.118pt}}
\multiput(346.00,302.17)(12.630,10.000){2}{\rule{0.330pt}{0.400pt}}
\multiput(360.00,313.58)(0.652,0.491){17}{\rule{0.620pt}{0.118pt}}
\multiput(360.00,312.17)(11.713,10.000){2}{\rule{0.310pt}{0.400pt}}
\multiput(373.00,323.58)(0.652,0.491){17}{\rule{0.620pt}{0.118pt}}
\multiput(373.00,322.17)(11.713,10.000){2}{\rule{0.310pt}{0.400pt}}
\multiput(386.00,333.59)(0.786,0.489){15}{\rule{0.722pt}{0.118pt}}
\multiput(386.00,332.17)(12.501,9.000){2}{\rule{0.361pt}{0.400pt}}
\multiput(400.00,342.59)(0.728,0.489){15}{\rule{0.678pt}{0.118pt}}
\multiput(400.00,341.17)(11.593,9.000){2}{\rule{0.339pt}{0.400pt}}
\multiput(413.00,351.59)(0.728,0.489){15}{\rule{0.678pt}{0.118pt}}
\multiput(413.00,350.17)(11.593,9.000){2}{\rule{0.339pt}{0.400pt}}
\multiput(426.00,360.59)(0.890,0.488){13}{\rule{0.800pt}{0.117pt}}
\multiput(426.00,359.17)(12.340,8.000){2}{\rule{0.400pt}{0.400pt}}
\multiput(440.00,368.59)(0.824,0.488){13}{\rule{0.750pt}{0.117pt}}
\multiput(440.00,367.17)(11.443,8.000){2}{\rule{0.375pt}{0.400pt}}
\multiput(453.00,376.59)(0.950,0.485){11}{\rule{0.843pt}{0.117pt}}
\multiput(453.00,375.17)(11.251,7.000){2}{\rule{0.421pt}{0.400pt}}
\multiput(466.00,383.59)(0.890,0.488){13}{\rule{0.800pt}{0.117pt}}
\multiput(466.00,382.17)(12.340,8.000){2}{\rule{0.400pt}{0.400pt}}
\multiput(480.00,391.59)(0.950,0.485){11}{\rule{0.843pt}{0.117pt}}
\multiput(480.00,390.17)(11.251,7.000){2}{\rule{0.421pt}{0.400pt}}
\multiput(493.00,398.59)(1.123,0.482){9}{\rule{0.967pt}{0.116pt}}
\multiput(493.00,397.17)(10.994,6.000){2}{\rule{0.483pt}{0.400pt}}
\multiput(506.00,404.59)(1.026,0.485){11}{\rule{0.900pt}{0.117pt}}
\multiput(506.00,403.17)(12.132,7.000){2}{\rule{0.450pt}{0.400pt}}
\multiput(520.00,411.59)(1.123,0.482){9}{\rule{0.967pt}{0.116pt}}
\multiput(520.00,410.17)(10.994,6.000){2}{\rule{0.483pt}{0.400pt}}
\multiput(533.00,417.59)(1.378,0.477){7}{\rule{1.140pt}{0.115pt}}
\multiput(533.00,416.17)(10.634,5.000){2}{\rule{0.570pt}{0.400pt}}
\multiput(546.00,422.59)(1.214,0.482){9}{\rule{1.033pt}{0.116pt}}
\multiput(546.00,421.17)(11.855,6.000){2}{\rule{0.517pt}{0.400pt}}
\multiput(560.00,428.59)(1.378,0.477){7}{\rule{1.140pt}{0.115pt}}
\multiput(560.00,427.17)(10.634,5.000){2}{\rule{0.570pt}{0.400pt}}
\multiput(573.00,433.59)(1.378,0.477){7}{\rule{1.140pt}{0.115pt}}
\multiput(573.00,432.17)(10.634,5.000){2}{\rule{0.570pt}{0.400pt}}
\multiput(586.00,438.60)(1.943,0.468){5}{\rule{1.500pt}{0.113pt}}
\multiput(586.00,437.17)(10.887,4.000){2}{\rule{0.750pt}{0.400pt}}
\multiput(600.00,442.60)(1.797,0.468){5}{\rule{1.400pt}{0.113pt}}
\multiput(600.00,441.17)(10.094,4.000){2}{\rule{0.700pt}{0.400pt}}
\multiput(613.00,446.60)(1.797,0.468){5}{\rule{1.400pt}{0.113pt}}
\multiput(613.00,445.17)(10.094,4.000){2}{\rule{0.700pt}{0.400pt}}
\multiput(626.00,450.61)(2.918,0.447){3}{\rule{1.967pt}{0.108pt}}
\multiput(626.00,449.17)(9.918,3.000){2}{\rule{0.983pt}{0.400pt}}
\multiput(640.00,453.60)(1.797,0.468){5}{\rule{1.400pt}{0.113pt}}
\multiput(640.00,452.17)(10.094,4.000){2}{\rule{0.700pt}{0.400pt}}
\multiput(653.00,457.61)(2.695,0.447){3}{\rule{1.833pt}{0.108pt}}
\multiput(653.00,456.17)(9.195,3.000){2}{\rule{0.917pt}{0.400pt}}
\put(666,460.17){\rule{2.900pt}{0.400pt}}
\multiput(666.00,459.17)(7.981,2.000){2}{\rule{1.450pt}{0.400pt}}
\put(680,462.17){\rule{2.700pt}{0.400pt}}
\multiput(680.00,461.17)(7.396,2.000){2}{\rule{1.350pt}{0.400pt}}
\put(693,464.17){\rule{2.700pt}{0.400pt}}
\multiput(693.00,463.17)(7.396,2.000){2}{\rule{1.350pt}{0.400pt}}
\put(706,466.17){\rule{2.900pt}{0.400pt}}
\multiput(706.00,465.17)(7.981,2.000){2}{\rule{1.450pt}{0.400pt}}
\put(720,467.67){\rule{3.132pt}{0.400pt}}
\multiput(720.00,467.17)(6.500,1.000){2}{\rule{1.566pt}{0.400pt}}
\put(733,468.67){\rule{3.132pt}{0.400pt}}
\multiput(733.00,468.17)(6.500,1.000){2}{\rule{1.566pt}{0.400pt}}
\put(746,469.67){\rule{3.373pt}{0.400pt}}
\multiput(746.00,469.17)(7.000,1.000){2}{\rule{1.686pt}{0.400pt}}
\put(799,469.67){\rule{3.373pt}{0.400pt}}
\multiput(799.00,470.17)(7.000,-1.000){2}{\rule{1.686pt}{0.400pt}}
\put(813,468.67){\rule{3.132pt}{0.400pt}}
\multiput(813.00,469.17)(6.500,-1.000){2}{\rule{1.566pt}{0.400pt}}
\put(826,467.67){\rule{3.132pt}{0.400pt}}
\multiput(826.00,468.17)(6.500,-1.000){2}{\rule{1.566pt}{0.400pt}}
\put(839,466.17){\rule{2.900pt}{0.400pt}}
\multiput(839.00,467.17)(7.981,-2.000){2}{\rule{1.450pt}{0.400pt}}
\put(853,464.17){\rule{2.700pt}{0.400pt}}
\multiput(853.00,465.17)(7.396,-2.000){2}{\rule{1.350pt}{0.400pt}}
\put(866,462.17){\rule{2.700pt}{0.400pt}}
\multiput(866.00,463.17)(7.396,-2.000){2}{\rule{1.350pt}{0.400pt}}
\put(879,460.17){\rule{2.900pt}{0.400pt}}
\multiput(879.00,461.17)(7.981,-2.000){2}{\rule{1.450pt}{0.400pt}}
\multiput(893.00,458.95)(2.695,-0.447){3}{\rule{1.833pt}{0.108pt}}
\multiput(893.00,459.17)(9.195,-3.000){2}{\rule{0.917pt}{0.400pt}}
\multiput(906.00,455.94)(1.797,-0.468){5}{\rule{1.400pt}{0.113pt}}
\multiput(906.00,456.17)(10.094,-4.000){2}{\rule{0.700pt}{0.400pt}}
\multiput(919.00,451.95)(2.918,-0.447){3}{\rule{1.967pt}{0.108pt}}
\multiput(919.00,452.17)(9.918,-3.000){2}{\rule{0.983pt}{0.400pt}}
\multiput(933.00,448.94)(1.797,-0.468){5}{\rule{1.400pt}{0.113pt}}
\multiput(933.00,449.17)(10.094,-4.000){2}{\rule{0.700pt}{0.400pt}}
\multiput(946.00,444.94)(1.797,-0.468){5}{\rule{1.400pt}{0.113pt}}
\multiput(946.00,445.17)(10.094,-4.000){2}{\rule{0.700pt}{0.400pt}}
\multiput(959.00,440.94)(1.943,-0.468){5}{\rule{1.500pt}{0.113pt}}
\multiput(959.00,441.17)(10.887,-4.000){2}{\rule{0.750pt}{0.400pt}}
\multiput(973.00,436.93)(1.378,-0.477){7}{\rule{1.140pt}{0.115pt}}
\multiput(973.00,437.17)(10.634,-5.000){2}{\rule{0.570pt}{0.400pt}}
\multiput(986.00,431.93)(1.378,-0.477){7}{\rule{1.140pt}{0.115pt}}
\multiput(986.00,432.17)(10.634,-5.000){2}{\rule{0.570pt}{0.400pt}}
\multiput(999.00,426.93)(1.214,-0.482){9}{\rule{1.033pt}{0.116pt}}
\multiput(999.00,427.17)(11.855,-6.000){2}{\rule{0.517pt}{0.400pt}}
\multiput(1013.00,420.93)(1.378,-0.477){7}{\rule{1.140pt}{0.115pt}}
\multiput(1013.00,421.17)(10.634,-5.000){2}{\rule{0.570pt}{0.400pt}}
\multiput(1026.00,415.93)(1.123,-0.482){9}{\rule{0.967pt}{0.116pt}}
\multiput(1026.00,416.17)(10.994,-6.000){2}{\rule{0.483pt}{0.400pt}}
\multiput(1039.00,409.93)(1.026,-0.485){11}{\rule{0.900pt}{0.117pt}}
\multiput(1039.00,410.17)(12.132,-7.000){2}{\rule{0.450pt}{0.400pt}}
\multiput(1053.00,402.93)(1.123,-0.482){9}{\rule{0.967pt}{0.116pt}}
\multiput(1053.00,403.17)(10.994,-6.000){2}{\rule{0.483pt}{0.400pt}}
\multiput(1066.00,396.93)(0.950,-0.485){11}{\rule{0.843pt}{0.117pt}}
\multiput(1066.00,397.17)(11.251,-7.000){2}{\rule{0.421pt}{0.400pt}}
\multiput(1079.00,389.93)(0.890,-0.488){13}{\rule{0.800pt}{0.117pt}}
\multiput(1079.00,390.17)(12.340,-8.000){2}{\rule{0.400pt}{0.400pt}}
\multiput(1093.00,381.93)(0.950,-0.485){11}{\rule{0.843pt}{0.117pt}}
\multiput(1093.00,382.17)(11.251,-7.000){2}{\rule{0.421pt}{0.400pt}}
\multiput(1106.00,374.93)(0.824,-0.488){13}{\rule{0.750pt}{0.117pt}}
\multiput(1106.00,375.17)(11.443,-8.000){2}{\rule{0.375pt}{0.400pt}}
\multiput(1119.00,366.93)(0.890,-0.488){13}{\rule{0.800pt}{0.117pt}}
\multiput(1119.00,367.17)(12.340,-8.000){2}{\rule{0.400pt}{0.400pt}}
\multiput(1133.00,358.93)(0.728,-0.489){15}{\rule{0.678pt}{0.118pt}}
\multiput(1133.00,359.17)(11.593,-9.000){2}{\rule{0.339pt}{0.400pt}}
\multiput(1146.00,349.93)(0.728,-0.489){15}{\rule{0.678pt}{0.118pt}}
\multiput(1146.00,350.17)(11.593,-9.000){2}{\rule{0.339pt}{0.400pt}}
\multiput(1159.00,340.93)(0.786,-0.489){15}{\rule{0.722pt}{0.118pt}}
\multiput(1159.00,341.17)(12.501,-9.000){2}{\rule{0.361pt}{0.400pt}}
\multiput(1173.00,331.92)(0.652,-0.491){17}{\rule{0.620pt}{0.118pt}}
\multiput(1173.00,332.17)(11.713,-10.000){2}{\rule{0.310pt}{0.400pt}}
\multiput(1186.00,321.92)(0.652,-0.491){17}{\rule{0.620pt}{0.118pt}}
\multiput(1186.00,322.17)(11.713,-10.000){2}{\rule{0.310pt}{0.400pt}}
\multiput(1199.00,311.92)(0.704,-0.491){17}{\rule{0.660pt}{0.118pt}}
\multiput(1199.00,312.17)(12.630,-10.000){2}{\rule{0.330pt}{0.400pt}}
\multiput(1213.00,301.92)(0.652,-0.491){17}{\rule{0.620pt}{0.118pt}}
\multiput(1213.00,302.17)(11.713,-10.000){2}{\rule{0.310pt}{0.400pt}}
\multiput(1226.00,291.92)(0.590,-0.492){19}{\rule{0.573pt}{0.118pt}}
\multiput(1226.00,292.17)(11.811,-11.000){2}{\rule{0.286pt}{0.400pt}}
\multiput(1239.00,280.92)(0.590,-0.492){19}{\rule{0.573pt}{0.118pt}}
\multiput(1239.00,281.17)(11.811,-11.000){2}{\rule{0.286pt}{0.400pt}}
\multiput(1252.00,269.92)(0.582,-0.492){21}{\rule{0.567pt}{0.119pt}}
\multiput(1252.00,270.17)(12.824,-12.000){2}{\rule{0.283pt}{0.400pt}}
\multiput(1266.00,257.92)(0.590,-0.492){19}{\rule{0.573pt}{0.118pt}}
\multiput(1266.00,258.17)(11.811,-11.000){2}{\rule{0.286pt}{0.400pt}}
\multiput(1279.00,246.92)(0.539,-0.492){21}{\rule{0.533pt}{0.119pt}}
\multiput(1279.00,247.17)(11.893,-12.000){2}{\rule{0.267pt}{0.400pt}}
\multiput(1292.00,234.92)(0.536,-0.493){23}{\rule{0.531pt}{0.119pt}}
\multiput(1292.00,235.17)(12.898,-13.000){2}{\rule{0.265pt}{0.400pt}}
\multiput(1306.00,221.92)(0.539,-0.492){21}{\rule{0.533pt}{0.119pt}}
\multiput(1306.00,222.17)(11.893,-12.000){2}{\rule{0.267pt}{0.400pt}}
\multiput(1319.00,209.92)(0.497,-0.493){23}{\rule{0.500pt}{0.119pt}}
\multiput(1319.00,210.17)(11.962,-13.000){2}{\rule{0.250pt}{0.400pt}}
\multiput(1332.00,196.92)(0.497,-0.494){25}{\rule{0.500pt}{0.119pt}}
\multiput(1332.00,197.17)(12.962,-14.000){2}{\rule{0.250pt}{0.400pt}}
\multiput(1346.00,182.92)(0.497,-0.493){23}{\rule{0.500pt}{0.119pt}}
\multiput(1346.00,183.17)(11.962,-13.000){2}{\rule{0.250pt}{0.400pt}}
\multiput(1359.58,168.80)(0.493,-0.536){23}{\rule{0.119pt}{0.531pt}}
\multiput(1358.17,169.90)(13.000,-12.898){2}{\rule{0.400pt}{0.265pt}}
\multiput(1372.58,154.81)(0.494,-0.534){25}{\rule{0.119pt}{0.529pt}}
\multiput(1371.17,155.90)(14.000,-13.903){2}{\rule{0.400pt}{0.264pt}}
\multiput(1386.58,139.80)(0.493,-0.536){23}{\rule{0.119pt}{0.531pt}}
\multiput(1385.17,140.90)(13.000,-12.898){2}{\rule{0.400pt}{0.265pt}}
\multiput(1399.58,125.67)(0.493,-0.576){23}{\rule{0.119pt}{0.562pt}}
\multiput(1398.17,126.83)(13.000,-13.834){2}{\rule{0.400pt}{0.281pt}}
\multiput(1412.58,110.81)(0.494,-0.534){25}{\rule{0.119pt}{0.529pt}}
\multiput(1411.17,111.90)(14.000,-13.903){2}{\rule{0.400pt}{0.264pt}}
\multiput(1426.58,95.54)(0.493,-0.616){23}{\rule{0.119pt}{0.592pt}}
\multiput(1425.17,96.77)(13.000,-14.771){2}{\rule{0.400pt}{0.296pt}}
\put(760.0,471.0){\rule[-0.200pt]{9.395pt}{0.400pt}}
\end{picture}

\caption{Rappresentazione grafica delle probabilit\`a descritte nel 
testo.}
\label{fig1}
\end{figure}

Poich\'e $P_{err2}$ \`e sempre maggiore o uguale a $P_{err1}$, e
$P_{err2}=P_{err1}$ solo con $p=0$~o~1 (vedi fig.~\ref{fig1}), se Sara
conta due volte i suoi soldi ha una probabilit\`a di
sbagliare\footnote{Cio\`e di sbagliare almeno una volta, ma se ottiene due
valori diversi, magari uno corretto e l'altro sbagliato o anche entrambi
sbagliati, Sara non pu\`o decidere subito quale sia quello corretto e si
trova costretta a contare almeno una terza volta.} che \`e maggiore della
probabilit\`a di sbagliare il valore della cifra totale se li contasse una
sola volta!

Nel caso in cui Sara ottenga due valori identici nelle due successive
o\-pe\-ra\-zio\-ni, la probabilit\`a di aver sbagliato in entrambi i casi
\`e $p^2$, mentre la probabilit\`a di avere in mano i due valori corretti
\`e $(1-p)^2$. Ovviamente Sara non pu\`o sapere con certezza in quale
delle due situazioni si trova.

Se invece nelle due o\-pe\-ra\-zio\-ni ottiene due valori diversi,
sicuramente almeno uno di questi \`e sbagliato. Ad esempio, la
probabilit\`a che in due conteggi Sara ottenga un valore corretto e uno
sbagliato \`e $2p(1-p)$. E questa \`e maggiore di $p$ se $p<\frac{1}{2}$
(vedi fig.~\ref{fig1}). Il che significa che pi\`u \`e piccola la
probabilit\`a che Sara ha di sbagliare, con pi\`u probabilit\`a si deve
preparare a contare i suoi soldi una terza volta.

In sostanza, se $p<\frac{1}{2}$ e Sara decide di contare i suoi soldi due
volte, allora si dovr\`a preparare a contarli una terza volta con una
probabilit\`a maggiore rispetto a quella di sbagliare il conteggio se li
contasse una volta sola!

La figura~\ref{fig1} riassume graficamente quanto \`e stato esposto sino
ad ora. Va comunque fatto notare che la curva che da la probabilit\`a di
ottenere due valori diversi in due conteggi successivi segue un andamento
diverso da quello indicato dalla funzione $2p(1-p)$, soprattutto per
$p>\frac{1}{2}$. Infatti, la funzione $2p(1-p)$ non ci da {\em tutta} la
probabilit\`a di ottenere due valori diversi in due conteggi successivi,
ma solo la probabilit\`a di ottenere un valore corretto e un valore
sbagliato. Se $p$ comincia ad essere significativamente grande (cio\`e
vicina ad 1), la probabilit\`a $2p(1-p)$ va a zero, ma non va certo a zero
la probabilit\`a di ottenere due valori diversi, entrambi sbagliati.
Quindi bisogna aspettarsi che la probabilit\`a di dover contare i soldi
una terza volta non vada a zero per $p$ che va a 1 (al contrario, essa
tende sicuramente ad 1) e, inoltre, che essa sia maggiore di $p$ anche per
valori pi\`u grandi di $p=\frac{1}{2}$. Tuttavia, fornire una stima di
quest'ultima probabilit\`a risulta essere un impresa non da poco, visto
che bisogna conoscere la forma esplicita di $p$ in funzione del numero
totale da contare e dall'entit\`a dell'errore commesso (differenza fra il
valore risultante dal conteggio e il valore esatto).

\subsection{Digressione filosofica}

Dati $M$ e $N$, cio\`e il numero di oggetti da contare e il numero di
volte che vengono contati, $p(\Delta_i)$ \`e definita come la
probabilit\`a di ottenere in un singolo conteggio un valore $E$ che
differisce da $M$ della quantit\`a $\Delta_i$, dove $\Delta_i=E-M$.
L'intervallo di valori (interi) che $\Delta_i$ pu\`o coprire \`e
$[\Delta_{min};\Delta_{max}]$, con $\Delta_{min}$ e $\Delta_{max}$
dipendenti in qualche modo anche da $M$ (cio\`e se $M$ \`e molto grande
\`e pi\`u probabile commettere un errore maggiore; nel contare pi\`u
oggetti \`e plausibilmente pi\`u difficile mantenere una certa
concentrazione). Notare che in linea di principio $\Delta_{min}$ \`e un
numero intero negativo e $\Delta_{max}$ positivo.  Quindi abbiamo che
$p(0)=1-p$, dove $p$ \`e la probabilit\`a di sbagliare introdotta nel
paragrafo precedente.

Ora, per ogni $N$ e $M$ esiste una probabilit\`a non nulla $P_{un}$ che
una persona che conta $N$ volte $M$ oggetti simili (monete, righe in una
lista di persone, schede elettorali in un'urna...) non sia in grado di
decidere quanti siano con precisione. E quindi potr\`a esistere sempre
qualcuno o un occasione in cui ognuno, per quante volte conti questi
oggetti, non sar\`a mai capace di dire con sicurezza quanti sono
realmente. Questo perch\'e chi conta non ottiene sempre lo stesso numero e
i diversi valori ottenuti (sia corretti che sbagliati) compaiono in gruppi
circa ugualmente numerosi. Non si \`e quindi in grado di decidere.

Tenendo conto di tutte le possibili combinazioni, $P_{un}$ assume la
seguente forma generale

\begin{equation}
P_{un}=\sum^{\Delta_{max}-\Delta_{min}+1}_{l\geq 2,\,\,\, k_i\simeq
k_j} \frac{N!}{\prod^l_{i=1}
k_i!}\sum_{\Delta_{max}\geq\Delta_{l}>\Delta_{l-1} > \dots > \Delta_{3} >
\Delta_{2}\geq \Delta_{min}} \prod^l_{i=1}p^{k_i}(\Delta_i), 
\end{equation}
con $\sum^l_{i=1} k_i = N$, dove $\Delta_{max}-\Delta_{min}+1$ \`e il
numero totale di valori diversi che possono essere ottenuti in un singolo
conteggio (quindi \`e incluso anche il caso in cui il valore ottenuto \`e
quello giusto, ovvero $\Delta_i=0$) e la prima sommatoria \`e eseguita in 
maniera tale che in $N$ conteggi consecutivi i diversi valori ottenuti 
compaiano in numero circa uguale, $k_i\simeq k_j$.

La seconda sommatoria \`e eseguita per {\em tutti} i valori di $\Delta_l$
che soddisfano la condizione $\Delta_{max}\geq\Delta_{l}>\Delta_{l-1} >
\dots > \Delta_{3} > \Delta_{2}\geq \Delta_{min}$.

Anche qui, per ottenere una stima di $P_{un}$ bisognerebbe conoscere
$p(\Delta_i)$ e, come abbiamo gi\`a detto prima, non \`e un'impresa
semplice.

Comunque, per ogni $N$ esister\`a un'occasione, per quanto improbabile, in
cui la quantit\`a da contare risulta inconoscibile all'essere umano che
tenta di contarla. Inoltre, un $M$ molto grande limita spesso la
possibilit\`a fisica di ripetere un conteggio e quindi impone un valore
molto piccolo per $N$: questo presumibilmente rende significativa la
probabilit\`a $P_{un}$.

In pi\`u, sarebbe interessante sapere se esistono $M$ e $N$ tali che,
anche per $p(\Delta_i)\simeq 0$, la probabilit\`a $P_{un}$ risulta essere
significativamente maggiore di zero.

\section{\bf Il paradosso del genio}

Si racconta\footnote{Si veda, ad esempio, W.E.~Deming, ``Out of the
crisis'' (MIT, Cambridge, 1986).} che un giorno, durante le fasi del
Progetto Manhattan (la co\-stru\-zio\-ne della bomba atomica), Enrico
Fermi chiese al Generale Robert Groves, capo del progetto, quale fosse la
definizione di `grande' generale. Groves rispose che ogni generale che
avesse vinto cinque battaglie consecutive poteva definirsi a buon diritto
`grande'. Fermi allora chiese quanti fossero i grandi generali della
storia e Groves rispose che ce n'erano circa tre ogni cento.

Fermi, facendo la semplicistica assunzione che nella maggior parte delle
battaglie le forze in campo tra le parti si equivalessero, assegn\`o alla
probabilit\`a di vincere una battaglia il valore di 1/2, cio\`e il $50\%$,
e calcol\`o che la probabilit\`a di vincere cinque battaglie
consecutive\footnote{Che \`e anche la probabilit\`a di ottenere 5 teste o
5 croci consecutive in 5 lanci di una moneta non viziata.} era pari a
$(1/2)^5=1/32$, cio\`e circa 3/100. ``Quindi ha ragione Generale''
replic\`o Fermi ``circa tre ogni cento. Probabilit\`a matematica, non
genio''.

L'aneddoto appena descritto non pu\`o non impressionare, visto che il
numero di `grandi' generali nella storia, secondo l'argomento di Fermi,
sarebbe lo stesso anche se si decidesse la vittoria delle battaglie
lanciando una monetina, piuttosto che combattendo!

Tuttavia, ci\`o finisce per suscitare una ulteriore riflessione. Se 3 su
100 \`e anche la probabilit\`a matematica nell'ipotesi di vittoria {\em
casuale}, quale avrebbe dovuto essere il valore della frazione di `grandi'
generali nella storia affinch\'e si fosse potuto parlare di esistenza del
`genio'?

La frequenza reale di `grandi' generali nella storia non pu\`o essere
inferiore a quella che si avrebbe per puro caso (1/32 appunto, o $\sim
3/100$). E questo \`e ovvio poich\'e il puro caso rappresenta una sorta di
limite inferiore; se i generali stessero con le mani in mano e non si
sforzassero di vincere, il risultato sarebbe comunque di 1/32. Allora,
l'unica possibilit\`a che questa frequenza ha di essere diversa dal valore
casuale \`e che sia superiore ad esso. Ma quanto maggiore? Pi\`u grande
\`e e meglio \`e?

E qui si pone un'altro dilemma:  maggiore \`e la frazione di `grandi'
generali fra tutti i generali, minore \`e il valore del termine `grande',
perch\'e in qualche senso risulta pi\`u facile essere `grande'. Ma allora
che significa essere `grande', essere un `genio' dopo queste premesse?

Pi\`u in generale, supponiamo che sia possibile associare una
probabilit\`a {\em a priori}, $p$, ad un evento `fuori dal comune' (ad
esempio, vittoria di 10 battaglie consecutive o di 2 premi Nobel da parte
di uno stesso individuo, scoperta di 2 teorie fondamentali della fisica in
campi diversi, etc.) partendo dall'ipotesi nulla, cio\`e che tutto avvenga
per caso, e opportunamente elaborata, dal punto di vista matematico, per
tenere conto della complessit\`a del problema specifico\footnote{Il fatto
che poi sia spesso praticamente impossibile calcolare $p$ non cambia il
punto del nostro discorso.}. Allora la frequenza media di esseri umani
`grandi' o `geni' \`e $f_0=Np$, dove $N$ \`e il numero che tiene conto di
tutti gli esseri umani vissuti in certo periodo storico. Ovviamente, $f_0$
\`e la frequenza nell'ipotesi casuale.

Ora, la frequenza reale di tali individui nella storia, $f_a$, non pu\`o
essere inferiore a $f_0$, per lo stesso motivo spiegato pi\`u sopra.
Inoltre, se $f_a=f_0$ si pu\`o dire anche qui, con Fermi, che si tratta di
probabilit\`a matematica, non di genio.

Allora, viste le premesse, affinch\'e si possa parlare di `genio' deve
essere necessariamente che $f_a>f_0$, cio\`e la frequenza reale di persone
che hanno raggiunto dei risultati eccezionali deve essere strettamente
maggiore di quella nell'ipotesi casuale. Ma `maggiore' \`e un concetto
vago. Pi\`u grande \`e e meglio \`e?

Anche in questo caso, maggiore \`e $f_a$, maggiore risulta essere la
frazione di persone `grandi' o `geni' nella popolazione totale
considerata e quindi minore risulta essere il `valore' della loro
grandezza, appunto perch\'e pi\`u comune in un senso preciso. Si pu\`o
applicare qui una sorta di Teoria dell'Utili\`a Marginale della
`grandezza' o `genialit\`a'. Nel campo dell'Economia la Teoria
dell'Utilit\`a Marginale afferma che l'utilit\`a di una porzione
aggiuntiva di bene materiale \`e legata alla quantit\`a totale di tale
bene gi\`a disponibile:  maggiore \`e tale quantit\`a, minore \`e
l'utilit\`a marginale della porzione aggiuntiva.

Potremmo dire, alternativamente, che se a raggiungere un risultato di
eccellenza fosse una frazione $f_{a_2}$ di persone minore di $f_{a_1}$,
noi percepiremmo il raggiungimento di tale risultato come pi\`u
eccezionale; ma questo processo potrebbe procedere a cascata fino al suo
limite inferiore, $f_0$ appunto, cio\`e
$f_{a_1}>f_{a_2}>\dots>f_{a_n}=f_0$.  Ci troviamo quindi di fronte ad un
paradosso.

Si sarebbe tentati di dire che l'idea comune di `genio', cio\`e di
spirito libero e creatore, novello Prometeo, che sfida la Natura e si
contrappone ad essa, sia in realt\`a uno di quei concetti vaghi creati
dall'uomo e che non abbia nessun riscontro oggettivo nella realt\`a.  
Tenendo conto di quanto \`e stato delineato qui, sembra piuttosto che il
`genio' sia l'essere umano che ha avuto l'occasione di trovarsi a vivere
in determinate circostanze e una determinata vita e questa possibilit\`a
\`e gi\`a meccanicamente prevista dalle leggi del caso: anche i `geni'
non sfuggono ad esso.

\section{\bf Democrazia della privazione}

Consideriamo un semplice modello di societ\`a costituita da $2N$
cittadini.  Supponiamo inoltre che questi $2N$ cittadini percepiscano uno
stipendio identico, pari a $s$ Euro.

Ora, cosa succederebbe se si mettesse al voto democratico\footnote{Cio\`e,
ricordiamo, tale che una qualsiasi riforma venga approvata se raggiunge
almeno il $50\% + 1$ dei voti.} la seguente riforma? Il testo della
riforma recita:

\begin{itemize}
\item una quantit\`a di denaro pari a $m$ viene tolta dallo stipendio di   
$N-k$ persone, con $k\geq 1$,

\item allo stipendio delle restanti $N+k$ persone viene aggiunta una
quantit\`a di denaro pari a 

$$\frac{(N-k)m}{N+k},$$ 
cio\`e il denaro tolto ai $N-k$ cittadini viene equamente ridistributo tra
i restanti $N+k$ cittadini.

\end{itemize}

Ora, poich\'e $N+k\geq N +1$, l'eventuale voto positivo del gruppo dei
$N+k$ cittadini corrisponde ad almeno il $50\% + 1$ dei voti e la riforma
\`e destinata a passare con certezza, nuocendo inevitabilmente alle altre
$N-k$ persone. Inoltre, se tu elettore non sai di quale gruppo andrai a
far parte, i fortunati $N+k$ o gli sfortunati $N-k$, e se si pu\`o
ragionevolmente assumere che la selezione di questi due gruppi sar\`a
casuale (o, almeno, si hanno motivi per ritenere che non ci siano
condizioni che favoriscano la tua immissione in una categoria piuttosto
che in un'altra) allora ti conviene votare {\em favorevolmente}, poich\'e
la tua probabilit\`a di essere tra gli $N+k$ fortunati \`e

$$\frac{N+k}{2N}= 0.5+\frac{k}{2N}>0.5,$$
cio\`e maggiore del $50\%$.

Per inciso, se questa modifica casuale di stipendio fosse un gioco a cui
si pu\`o scegliere di partecipare, questo sarebbe un gioco a guadagno
zero. In un gioco in cui si pu\`o vincere o perdere una somma di denaro
secondo una certa probabilit\`a, il guadagno $G$ \`e definito come la
somma della probabilit\`a di vincere moltiplicata per la somma vinta e
della probabilit\`a di perdere moltiplicata per la somma persa (la somma
persa deve comparire con un segno meno). \`E in sostanza una media del
guadagno. Per chi partecipa al gioco come cliente, esempi tipici di giochi
a guadagno negativo sono: il lotto, il superenlotto e i giochi da Casin\`o
in genere (in questi casi il guadagno medio \`e positivo per i gestori dei
giochi). Nel nostro caso, il guadagno medio di ogni persona \`e cos\`{\i}
calcolato

$$G=\underbrace{\frac{N+k}{2N}}_{P_{\textrm{\tiny vittoria}}}
\underbrace{\frac{(N-k)m}{N+k}}_{\textrm{\tiny Somma 
vinta}}+\underbrace{\frac{N-k}{2N}}_{P_{\textrm{\tiny 
sconfitta}}}\underbrace{(-m)}_{\textrm{\tiny Somma persa}}=0,$$
e quindi rientra in quel gruppo di scommesse che la Teoria dei Giochi 
definisce {\it fair}, cio\`e oneste.

Per concludere, nell'ultima sezione si \`e voluto descrivere brevemente
questo esempio ideale per mostrare come possa esistere una linea razionale
e democratica di comportamento che ci favorisca a discapito degli altri.
Chi vota favorevolmente sapendo di essere fra i fortunati $N+k$ cittadini
esprime democraticamente il suo parere e non prevarica nessuno. In pi\`u,
abbiamo anche visto che se tu non sapessi di essere fra gli $N+k$ o fra
gli $N-k$, razionalmente ti converrebbe comunque votare in maniera
favorevole alla riforma.

Questo \`e un classico esempio (ed anche il pi\`u banale) di come ci\`o
che \`e perfettamente razionale per il singolo, ed \`e complessivamente
democratico, pu\`o essere nocivo per una porzione consistente della
collettivit\`a. E va da se che per impedire che si verifichino situazioni
del genere, nelle quali un comportamento razionale e democratico porta
paradossalmente a nuocere a molti altri, bisogna essere attenti a non
proporre riforme di questo tipo, o comunque a cercare di mantenere una
visione globale del problema.

\begin{center}
$\ast\ast\ast$
\end{center}
Gli autori ringraziano Fulvio~De~Cicco per aver letto una prima bozza
dell'articolo ed aver contribuito criticamente al suo miglioramento.

\end{document}